\documentclass[12pt]{amsart}

\usepackage{amssymb, amsmath, xcolor}
\usepackage{import}
\usepackage{tikz}
\usepackage{array}
\tikzset{filled/.style={minimum width=5pt,inner sep=0pt,circle,fill=black}}
\usepackage{subcaption}
\usepackage{mathtools}
\usepackage[utf8]{inputenc}
\usepackage{float}
\usepackage{longtable}
\usepackage{todonotes}
\usepackage[square,numbers]{natbib}
\newtheorem{theorem}{Theorem}[section]
\newtheorem{lemma}[theorem]{Lemma}
\newtheorem{corollary}[theorem]{Corollary}

\theoremstyle{definition}
\newtheorem{definition}[theorem]{Definition}

\theoremstyle{remark}
\newtheorem{remark}[theorem]{Remark}

\numberwithin{equation}{section}
\numberwithin{figure}{section}

\newcommand{\mpd}{\operatorname{mpd}}
\newcommand{\N}{\mathbb{N}}
\newcommand{\Z}{\mathbb{Z}}
\title[Criticality, Additivity, Join, Pseudoachromatic Number]{On Criticality and Additivity of the Pseudoachromatic Number Under Join}

\author[Meddaugh, Sepanski, Venkataraman]{
	Jonathan Meddaugh, Mark R. Sepanski, Yegnanarayanan Venkataraman }
\thanks{Communicating author: Y. Venkataraman}

\date{\today}

\address{
	Meddaugh and Sepanski:
	Department of Mathematics,
	Baylor University,
	Sid Richardson Building,
	1410 S.~4th Street,
	Waco, TX 76706, USA; \\
	Venkataraman: 
	Kalasalingam Academy of Research and Educatioon, Deemed to be University,
	Srivilliputhur, Tamil Nadu 626126, India}

\email{\texttt{Jonathan\_Meddaugh@baylor.edu, Mark\_Sepanski@baylor.edu, prof.yegna@gmail.com}}

\begin{document}
	
	

	\begin{abstract}
		A vertex coloring of a graph is said to be pseudocomplete if, for any two distinct colors, there exists at least one edge with those two colors as its end vertices.  The pseudoachromatic number of a graph is the greatest number of colors possible used in a pseudocomplete coloring. 
		This paper studies properties relating to additivity of the pseudoachromatic number under the join. Errors from the literuature are corrected and the notion of weakly critical is introduced in order to study the problem.
	\end{abstract}

	\maketitle
	\tableofcontents
	\section{Introduction}

	Some of the oldest problems in graph theory study colorings subject to various constraints. A coloring of the vertices of a graph is called \emph{pseudocomplete} if every pair of disjoint colors are adjacent via at least one edge. The maximum possible number of colors used in a pseudocomplete coloring of a graph $G$ is called the \emph{pseudoachromatic number}, $\Psi(G)$. Though not studied here, requiring the colorings be proper gives rise to the well studied achromatic number, \emph{e.g.}, \cite{HS1970}.
	
	The pseudoachromatic number was first used by Harary \emph{et al.}, \cite{HFS1967}, and Gupta, \cite{Gupta1968}.
	See \cite{Yegnanarayanan2000} for additional comments on the origin of the name.
	For some of the history of work on $\Psi$, see the following:
	\cite{Sampathkumar1976PartitionGA}, 
	\cite{B1979}, 
	\cite{Balakrishnan1998ExtremalGI}, 
	\cite{Yegnanarayanan2000OnTE}, 
	\cite{Edwards2000AchromaticNV}, 
	\cite{Yegnanarayanan2001GraphCA}, 
	\cite{Yegnanarayanan2002}, 
	\cite{Chen2008OnTP},
	\cite{Yegnanarayanan2009}, 
	\cite{Yegnanarayanan2013OnPA},
	\cite{LL2011},
	\cite{AGMJRCS2011},
	\cite{AGMJRCS2014},
	\cite{ACLR2016},
	\cite{AGMJRCS2017}, and
	\cite{AGMJRC2018}.
	
	Initial definitions and background results are found in Section \S\ref{sec: init defs}, including the definition of \emph{criticality}, Definition \ref{def: non-critical values and critical}, which plays a role in some aspects of $\Psi$ being additive under join. 
	
	Section \S\ref{sec: crit} examines additivity of $\Psi$ and preservation of criticality under the join.
	On this front, for graphs $G$ and $H$, there is an error in the literature that claims that $G$ and $H$ are critical if and only if
	\begin{equation} \label{eqn: intro psi additive}
		\Psi(G \,\nabla H) = \Psi(G) + \Psi(H).
	\end{equation}
	Though sufficient (Theorem \ref{thm: G join H psi additive when critical}), in fact, the converse is not true, see Remark \ref{rmk: erroneous claim on psi additve implies critical}. 
	Similarly (Theorem \ref{thm: crit join crit is crit}), when $G$ and $H$ are critical, so is $G\,\nabla H$, though the converse is also false, see Remark \ref{rmk: false: converse of crit join crit is crit}. 
	Finally apriori bounds and structural results are given in Theorems \ref{thm: upper lower bounds for general psi of join} and \ref{thm: crit lots of edges}.
	
	Section \S\ref{sec: weak crit} introduces the notion of \emph{weakly critical} in Definition \ref{def: non-critical values and critical}. 
	A graph invariant equivalence is given in Theorem \ref{thm: weak critical numerical equivalence} and a structural equivalence is given in Theorem \ref{thm: alt def of not weakly critical}.
	One of the main results of this paper, Theorem \ref{thm: psi add gives at least a weak crit on one}, shows that Equation \ref{eqn: intro psi additive} implies that either $G$ or $H$ is weakly critical.
	
	Section \S\ref{sec: nabla^k G} studies the above topics in the context of $\nabla^k G$. 
	For $k\geq 2$, Theorem \ref{thm: k-join of G critical when it can be} shows that $\nabla^k G$ is critical precisely when
	$k\,(\omega(G) + |V|)$ is even.
	More generally, Theorem \ref{thm: k-join of G weakly critical for the rest}
	shows that $\nabla^k G$ is always weakly critical.
	Finally, additivity or near-additivity of $\Psi$ on $\nabla^k G$ is determined in Theorems \ref{thm: crit if kpsi = psi nabla k with parity even condition} and \ref{thm: w crit if kpsi = psi nabla k with parity odd condition}.

	\section{Initial Definitions and Background} \label{sec: init defs}
	
	We write $\N$ for the nonnegative integers and $\Z^+$ for the positive integers. 
	In this paper, $G=(V,E)$ is a simple, finite graph. If multiple graphs may lead to ambiguity, we will write $G=(V_G,E_G)$. If the vertices of $G$ are equipped with a labeling by some $C\subseteq\Z$, we generally write $\ell:V \twoheadrightarrow C$ for the coloring.
	Finally, we write $\omega$ for the \emph{clique number} and $\nabla$ for the graph theoretic \emph{join}.
	
	We begin with one of the central definitions of the paper.
	
	\begin{definition} \label{def: pseudocomplete and Psi}
		A \emph{pseudocomplete coloring} of a graph $G$ is a coloring of the vertices of $G$, $\ell:V \rightarrow C$, so that, for any distinct $i,j\in C$, there exists $uv\in E$ so that 
		$\{i,j\} = \{\ell(u),\ell(v)\}$.
		Note that the coloring here need not be proper.
		
		The \emph{pseudoachromatic number} of a graph $G$, written $\Psi(G)$, is the greatest possible number of colors employed in a pseudocomplete coloring of $G$.
	\end{definition}
	
	If $K$ is a maximal clique of $G$, then coloring $K$ with distinct colors and the rest of $G$ arbitrarily gives a pseudocomplete coloring. In addition, the defintion of a pseudocomplete coloring requires at least $\binom{\Psi(G)}{2} \leq |E|$. As a result, 
	\[ \omega(G) \leq \Psi(G) \leq \frac{1+\sqrt{1+8|E|}}{2}.\]
	
	Moreover importantly, the following upper bound is known.
	\begin{lemma}[\cite{BRY2003}, Lemma 2] \label{lem: psi vs clique number inequality}
		Let $G$ be a graph. Then
		\[ \Psi(G) \leq \left\lfloor \frac{\omega(G) + |V|}{2} \right\rfloor. \]
	\end{lemma}
	
	In addition, there is a lower bound for $\Psi$ of the join of graphs.
	\begin{lemma}[\cite{BRY2003}, Corollary 4] \label{lem: psi of join geq sum of psi}
		Let $G$ and $H$ be graphs. Then
		\[ \Psi(G) + \Psi(H) \leq \Psi(G \,\nabla H). \]
	\end{lemma}
	
	Determining when the inequalities in Lemmas \ref{lem: psi vs clique number inequality} and \ref{lem: psi vs clique number inequality} become equalities leads to the following definitions. 
	
	\begin{definition} \label{def: mpd_G(k)}
		Let $G$ be a graph and $k\in\N$ with $k\leq |V|$. Define the \emph{minimal $psi$-drop function} as 
		\[ \mpd_G(k) = \min\{ \Psi(G) - \Psi(G\setminus X) \mid X\subseteq G, |X| = k \}. \]    
		
		If $X \subseteq G$ satifies $|X| = k$ and 
		$\mpd_G(k) = \Psi(G) - \Psi(G\setminus X)$, we call $X$ a \emph{realizing subgraph} for $\mpd_G(k)$.
	\end{definition}
	
	Note that 
	$ 0 \leq \mpd_G(k) \leq k $
	as $\Psi$ can drop by at most one when removing a vertex.
	
	The next definition is a reformulation of the one found in \cite{BRY2003}.

	\begin{definition} \label{def: non-critical values and critical}
		We say $G$ is \emph{critical} if 
		\[  \mpd_G(k) \geq \left\lceil \frac{k}{2} \right\rceil \]
		for all $0\leq k \leq |G|$.
%
%
	\end{definition}   
	
	There is a useful known equivalence for criticality in terms of maximizing the inequality in Lemma \ref{lem: psi vs clique number inequality} without the floor function.
	\begin{lemma}[\cite{BRY2003}, Lemma 10] \label{lem: crit if 2psi=o+n}
		Let $G$ be a graph. Then $G$ is critical if and only if 
		\[ \Psi(G) = \frac{\omega(G) + |V|}{2}. \]
	\end{lemma}

	\section{Criticality and Additivity of $\Psi$ Under Joins} \label{sec: crit}

	We begin with an apriori bound on $\Psi(G\,\nabla H)$
	between a minimal sum involving a clique number and vertex count and an average sum of clique numbers and vertex counts.
	\begin{theorem} \label{thm: upper lower bounds for general psi of join}
		Let $G$ and $H$ be graphs. Let 
		\[ m = \min\{\omega(G) + |V_H|,\,\, \omega(H) + |V_G|\}.\]
		Then
		\[ m
		\leq \Psi(G\,\nabla H) 
		\leq \left\lfloor
		\frac{\omega(G) + \omega(H)}{2} + \frac{|V_G| + |V_H|}{2}
		\right\rfloor. \]
	\end{theorem}
	
	\begin{proof}
		The upper bound comes from Lemma \ref{lem: psi vs clique number inequality} and the additivity of $\omega$ under join.
		
		For the lower bound, let $K_G$ and $K_H$ be maximal cliques of $G$ and $H$, respectively. Color $K_G$ and $K_H$ with $\omega(G) + \omega(H)$ distinct colors. Let $X_G = G\setminus K_G$ and $X_H = H\setminus K_H$. 
		After relabeling, we may assume that 
		\[ |V_G| - \omega(G) =|X_G| \leq |X_H| = |V_H| - \omega(H). \]
		Color both $X_G$ and $X_H$ with an additional identical $|V_G| - \omega(G)$ colors. It is straightforward to verify that this is a pseudocomplete coloring of $G\,\nabla H$ with 
		\[ (\omega(G) + \omega(H))+ (|V_G| - \omega(G)) = \omega(H) + |V_G|\]
		colors.
		After noting that $\omega(H) + |V_G| \leq \omega(G) + |V_H|$, we are done.
	\end{proof}
		
	Next we give a criteria for $\Psi$ to be additive over the join. 
	
	\begin{theorem} \label{thm: criteria for Psi to be additive over the join via mpd}
		Let $G, H$ be graphs. Then
		\[ \Psi(G\,\nabla H) = \Psi(G) + \Psi(H) \]
		if and only if
		\[ \mpd_G(k) + \mpd_H(k) \geq k \]
		for all $0\leq k \leq \min\{|G|, |H|\}$.
		
		Moreover, in that case, there exists a $k$, $1\leq k \leq \min\{|G|, |H|\}$, so that 
		\[\mpd_G(k) + \mpd_H(k) = k. \]
	\end{theorem}
	
	\begin{proof}
		Let $\ell$ be a pseudocomplete coloring of $G\,\nabla H$ with $\Psi(G\,\nabla H)$ colors. 
		Let $C = \ell(V_G) \cap \ell(V_H)$. By recoloring if necessary, we may assume that every color in $C$ appears exactly once in $G$ and once in $H$. 
		Define induced subgraphs of $G$ and $H$ by $X_G =\ell^{-1}(C) \cap V_G$ and $X_H =\ell^{-1}(C) \cap V_H$, respectively.
		We see that $|C| = |X_G| = |X_H|$. Note that $\ell$ restricted to $G\setminus X_G$ is pseudocomplete and satisfies $|\ell(G\setminus X_G)|=\Psi(G\setminus X_G)$ (and similar for $H\setminus X_H$).
		By construction, it follows that
		\[ \Psi(G\,\nabla H) = \Psi(G\setminus X_G) + |C| + \Psi(H\setminus X_H).\]
		
		From the observations that 
		$\Psi(G) \geq \Psi(G\setminus X_G) + \mpd_G(|C|)$ and
		$\Psi(H) \geq \Psi(H\setminus X_H) + \mpd_H(|C|)$, it follows that
		\[ \Psi(G) + \Psi(H) \geq \Psi(G\setminus X_G) + \Psi(H\setminus X_H) + \mpd_H(|C|) + \mpd_G(|C|)\]
		so that
		\[ \Psi(G) + \Psi(H) \geq  \Psi(G\,\nabla H) + (\mpd_H(|C|) + \mpd_G(|C|) -|C|). \]
		From this, we see that the condition $\mpd_G(k) + \mpd_H(k) \geq k$ for all $k$ forces 
		$\Psi(G) + \Psi(H) \geq  \Psi(G\,\nabla H)$,
		with equality for $k=|C|$,
		and, by Lemma \ref{lem: psi of join geq sum of psi}, additivity of $\Psi$ over the join.
		
		Conversely, suppose that $\mpd_G(k_0) + \mpd_H(k_0) < k_0$ for some $1\leq k_0 \leq \min\{|G|, |H|\}$.
		Choose realizing subgraphs for $\mpd_G(k_0)$ and $\mpd_H(k_0)$, 
		$X_G$ and $X_H$, respectively. 
		It follows 
		that  $\mpd_G(k_0) = \Psi(G) - \Psi(G\setminus X_G)$, $\mpd_H(k_0) = \Psi(H) - \Psi(H\setminus X_H)$, and  $|X_G| = |X_H| = k_0$. Then
		\[ \Psi(G) + \Psi(H) < \Psi(G\setminus X_G) + \Psi(H\setminus X_H) + k_0. \]
		
		We can pseudocompletely color $G\,\nabla H$ as follows: color $G\setminus X_G$ pseudocompletely with $\Psi(G\setminus X_G)$ colors, color $H\setminus X_H$ pseudocompletely with  $\Psi(H\setminus X_H)$ additional colors, and color both $X_G$ and $X_H$ with $k_0$ further colors. Then we have  \[ \Psi(G\setminus X_G) + \Psi(H\setminus X_H) + k_0 \leq  \Psi(G\,\nabla H)\] 
		and so
		\[ \Psi(G) + \Psi(H) <  \Psi(G\,\nabla H). \]
	\end{proof}

	From Theorem \ref{thm: criteria for Psi to be additive over the join via mpd}, we immediately get the following result.

	\begin{theorem}  \label{thm: G join H psi additive when critical}
		Let $G$ and $H$ be critical graphs. Then
		\[ \Psi(G \,\nabla H) = \Psi(G) + \Psi(H). \]
	\end{theorem}

	\begin{remark} \label{rmk: erroneous claim on psi additve implies critical}
		Though erroneously claimed in \cite{BRY2003} Corollary 6, the converse of Theorem \ref{thm: G join H psi additive when critical} is not true which invalidates some of the results and proofs in \cite{BRY2003}.
		The failure of the converse can be seen with $P_3$ and $C_8$. Using Lemma \ref{lem: crit if 2psi=o+n} and the straightforward facts that $\Psi(P_3) = 2$ and $\Psi(C_8) = 4$ 
		, neither $P_3$ nor $C_8$ is critical. However, a
		labeling of the join using $\{1,2,3\}$ and $\{1,4,5,6,4,4,4,4\}$ and Lemma \ref{lem: psi vs clique number inequality} show that $\Psi(P_3\,\nabla C_8) = 6$ so that 
		$\Psi(P_3\,\nabla C_8) \allowbreak =  \Psi(P_3) + \Psi(C_8)$.
		See Theorem \ref{thm: psi add gives at least a weak crit on one} below for as close to a converse as possible.
	\end{remark}

	Theorem \ref{thm: G join H psi additive when critical} allows us to show that criticality is preserved under join.
	\begin{theorem} \label{thm: crit join crit is crit}
		If $G$ and $H$ are critical graphs, then $G\,\nabla H$ is also critical.
	\end{theorem}

	\begin{proof}
		Use Lemma \ref{lem: crit if 2psi=o+n} and Theorem \ref{thm: G join H psi additive when critical} to see that $G\,\nabla H$ is critical if and only if
		\[2(\Psi(G) + \Psi(H))= (\omega(G)+\omega(H)) + (|V_G|+|W_H|)\]
		and we are done.
	\end{proof}

	\begin{remark} \label{rmk: false: converse of crit join crit is crit}
		The converse of Theorem \ref{thm: crit join crit is crit} is false. To see this, recall that $P_3$ is not critical. However, a labeling of $P_3\,\nabla P_3$ with $\{1,2,3\}$ and $\{1,4,5\}$ on each $P_3$ shows that $\Psi(P_3\,\nabla P_3) \geq 5$. Equality and the fact that
		$P_3\,\nabla P_3$ is critical then follow from 
		Lemma \ref{lem: psi vs clique number inequality} and 
		Lemma \ref{lem: crit if 2psi=o+n}.
		Note that this shows that the partial converse to Theorem \ref{thm: crit join crit is crit} is false even when restricting to the case of $G=H$. Of note, see the comment after Equation \ref{eqn: Psi of GjoinG} and Theorem \ref{thm: k-join of G critical when it can be}.
	\end{remark}

	We end this section with a general constraint on critical graphs and their pseudocomplete coloring.
	\begin{theorem} \label{thm: crit lots of edges}
		If $G$ is critical, then there is a pseudocomplete coloring of $G$ in which $\omega(G)$ colors are used exactly once on a maximal clique and $\frac{|V|-\omega}{2}$ colors are used exactly twice on the complement of the maximal clique. Moreover,
		\[ |E| \geq \frac{(|V|+\omega) \, (|V|+\omega -2)}{8}.\]
	\end{theorem}
	
	\begin{proof}
		Let $\ell$ be a pseudocomplete coloring of $G$ with the colors 
		$S=\{1,2,\ldots,\Psi(G)\}$. For $k\in\Z^+$, let
		\[ m_k = \{ i\in S \mid |\ell^{-1}(i)| = k\} \]
		and $n_k = |m_k|$. As $\ell$ is pseudocomplete, 
		$K = \bigcup_{i\in m_1} \ell^{-1}(i)$
		is a clique of size $n_1$.
		
		Moreover,
		\begin{align*}
			|V| &= \sum_{k\in\Z^+} k n_k \\
			&= n_1 + 2\sum_{k\geq 2} n_k + \sum_{k\geq 3} (k-2) n_k \\
			&= n_1 + 2(\Psi(G) - n_1) + \sum_{k\geq 3} (k-2) n_k \\
			&= 2\Psi(G) - n_1 + \sum_{k\geq 3} (k-2) n_k \\
			&\geq 2\Psi(G) - n_1 \\
			&\geq 2\Psi(G) - \omega(G).
		\end{align*}
		By Lemma \ref{lem: crit if 2psi=o+n}, $G$ is critical if and only if $|V| = 2\Psi(G) - \omega(G)$. From the equations above, we see $G$ is critical if and only if we have equality at each step. In particular, if and only if $m_k = \emptyset$ for $k\geq 3$ and $n_1 = \omega(G)$.
		
		As a result, $K$ is a maximal clique and its colors are used exactly once. In addition, the colors of $X = G\setminus K$ are used exactly twice. Therefore $X$ may be broken into two sets of equal size, $X_1$ and $X_2$, with the colors of $X_1$ used exactly once in $X_1$ and exactly once in $X_2$.
		
		Finally, write $\Tilde{X}$ for the contraction of $G$ that identifies each vertex in $X_1$ with the vertex in $X_2$ of the same color.
		As $\ell$ is a pseudocomplete coloring, we see that $\Tilde{X}$ is the complete graph on $\omega(G) + \frac{|V|-\omega(G)}{2}$ vertices and we are done.
	\end{proof}

	\section{Weak Criticality} \label{sec: weak crit}

	As demonstrated in Remark \ref{rmk: erroneous claim on psi additve implies critical}, criticality is not necessary for additivity of $\Psi$ over the join. Indeed, Theorem \ref{thm: criteria for Psi to be additive over the join via mpd} shows that additivity of $\Psi$ for a graph pair depends on on subtle structures in both graphs, rendering a simple characterization of additive pairs out of the question. Despite that, in this section we demonstrate that additivity of $\Psi$ requires at least one graph to have a weak form of criticality.

	The following is a subtle tweak of Definition \ref{def: non-critical values and critical}. Note the change from the ceiling function to the floor function.
	
	\begin{definition} \label{def: weakly critical}
		We say $G$ is \emph{weakly critical} if 
		\[  \mpd_G(k) \geq \left\lfloor \frac{k}{2} \right\rfloor \]
		for all $0\leq k \leq |G|$.
%
%
	\end{definition}   
	
	Note that Definitions \ref{def: non-critical values and critical} and \ref{def: weakly critical} imply that every critical graph is weakly critical. 	The converse to this is not true as $P_3$ is weakly critical but not critical.

	We immediately get the following by applying Theorem \ref{thm: criteria for Psi to be additive over the join via mpd} and noting that  $\left\lfloor \frac{k}{2} \right\rfloor + \left\lceil \frac{k}{2} \right\rceil=k$ for all $k$.
	
	\begin{corollary} \label{cor: WC plus C}
		Let $G$ and $H$ be graphs with $G$ critical and $H$  weakly critical. Then 
		\[\Psi(G\,\nabla H)=\Psi(G)+\Psi(H).\]
	\end{corollary}

	The following results concerning weak criticality are analogues of various results for criticality. The first is the analogue of Lemma \ref{lem: crit if 2psi=o+n} for being weakly critical and avoids the parity constraint of criticality.
	\begin{theorem} \label{thm: weak critical numerical equivalence}
		Let $G$ be a graph. Then $G$ is weakly critical if and only if 
		\[ \Psi(G) = \left\lfloor \frac{\omega(G) + |V|}{2} \right\rfloor. \]
	\end{theorem}
	
	\begin{proof}
		By Lemma \ref{lem: psi vs clique number inequality}, it suffices to prove that 
		\begin{equation} \label{eqn: inequality in weak crit numerical thm}
			\Psi(G) < \left\lfloor \frac{\omega(G) + |V|}{2} \right\rfloor
		\end{equation}
		happens if and only if G is not weakly critical.
		
		First suppose that Equation \ref{eqn: inequality in weak crit numerical thm} holds. Choose $K\subseteq G$ to be a maximal clique and let $X = G\setminus K$. Then
		\begin{align*}
			\Psi(G) &< \left\lfloor \frac{\omega(G) + |V|}{2} \right\rfloor \\
			&=  \omega(G) + \left\lfloor \frac{|V| - \omega(G)}{2} \right\rfloor \\
			&= \Psi(G\setminus X) + \left\lfloor \frac{|X|}{2} \right\rfloor.
		\end{align*}
		Therefore $\mpd_G(|X|)\leq \Psi(G)-\Psi(G\setminus X)<\left\lfloor \frac{|X|}{2} \right\rfloor$ and $G$ is not weakly critical. 
		
		Conversely, suppose there exists $G$ is not weakly critical. Then there exists $k\leq|G|$ with $\mpd_G(k)<\left\lfloor \frac{|X|}{2} \right\rfloor$. Let $X$ be a realizing subgraph for $\mpd_G(k)$. Then, using Lemma \ref{lem: psi vs clique number inequality},
		\begin{align*}
			\Psi(G) &< \Psi(G\setminus X) + \left\lfloor \frac{|X|}{2} \right\rfloor \\
			&\leq  \left\lfloor \frac{\omega(G\setminus X) + (|V| - |X|)}{2} \right\rfloor
			+ \left\lfloor \frac{|X|}{2} \right\rfloor \\
			&\leq \left\lfloor \frac{\omega(G) + |V| - |X|}{2} \right\rfloor
			+ \left\lfloor \frac{|X|}{2} \right\rfloor \\
			&\leq \left\lfloor \frac{\omega(G) + |V|}{2} \right\rfloor
		\end{align*}
		and we are done.
	\end{proof}

	From Theorem \ref{thm: weak critical numerical equivalence} and Lemma \ref{lem: crit if 2psi=o+n}, observe that if 
	$\omega(G) + |V|$ is even, then $G$ is critical if and only if it is weakly critical. 
	However, when $\omega(G) + |V|$ is odd, $G$ is not critical, but $G$ may still be weakly critical. The graph $P_3$ is such an example.
	
	We now give the analogue of Theorem \ref{thm: crit lots of edges} for weakly critical.
	\begin{theorem} \label{thm: weakly crit lots of edges}
		Let $G$ be weakly critical but not critical.
		Then there is a pseudocomplete coloring of $G$ so that either
		\begin{enumerate}
			\item $\omega(G)$ colors are used exactly once on a maximal clique, a single color is used exactly three times in the complement of the clique, and $\frac{|V|-\omega-3}{2}$ colors are used exactly twice on remaining vertices or
			\item$\omega(G)-1$ colors are used exactly once on a clique of size $\omega(G) -1$ and $\frac{|V|-\omega+1}{2}$ colors are used exactly twice on complement of the clique.
		\end{enumerate}
		Moreover,
		\[ |E| \geq \frac{(|V|+\omega-1) \, (|V|+\omega -3)}{8}.\]
	\end{theorem}
	
	\begin{proof}
		This result follows from the proof of Theorem \ref{thm: crit lots of edges}. As Theorem \ref{thm: weak critical numerical equivalence} shows that weakly critical in this setting is the same as requiring
		\[ |V| + \omega(G) - 2\Psi(G) -1 = 0,\]
		the proof of Theorem \ref{thm: crit lots of edges} shows that $n_k =0$ for $k\geq 4$. Moreover, either
		$n_1 = \omega$ and $m_3 = 1$ or $n_1 = \omega(G) -1$ and $m_3 =0$.
	\end{proof}
	
	The remainder of this section is devoted to demonstrating that weak criticality is required of at least one graph in a $\Psi$-additive graph pair.
	
	The following result will be an important technical tool.
	\begin{theorem} \label{thm: alt def of not weakly critical}
		Let $G$ be a graph. Then $G$ is not weakly critical if and only if 
		there exists induced subgraphs $M_1 \subseteq M_2 \subseteq G$ so that
		\begin{enumerate}
			\item $|V_{M_2 \setminus M_1}| = 2$,
			\item $\Psi(M_1) = \Psi(M_2)$,
			\item $\Psi(G) = \Psi(M_2) 
			+ \left\lfloor \frac{|V_{G\setminus M_2}|}{2} \right\rfloor$.
		\end{enumerate}
		Moreover, if $\ell$ is a maximal pseudocomplete coloring of $G$, 
		there exists $C\subseteq V_G$ with $|C| = \left\lfloor \frac{|V_{G\setminus M_1}|}{2} \right\rfloor + 1\geq 2$ and
		$\ell(G) = \ell(G\setminus C)$.
	\end{theorem}
	
	\begin{proof}
		Suppose $G$ is not weakly critical and, by Definition \ref{def: weakly critical}, choose $M_1$ to be a maximal subgraph of $G$ satisfying
		\[ \Psi(G) \leq \Psi(M_1) 
		+ \left\lfloor \frac{|V_{G\setminus M_1}|}{2} \right\rfloor -1.\]
		As $\Psi(G) \geq \Psi(M_1)$, $|V_{G\setminus M_1}| \geq 2$. Choose $M_2$ to be any induced subgraph of $G$ satisfying $M_1 \subseteq M_2$ with 
		$|V_{M_2 \setminus M_1}| = 2$.
		
		By maximality, we get
		\begin{align*}
			\Psi(G) &\geq \Psi(M_2) 
			+ \left\lfloor \frac{|V_{G\setminus M_2}|}{2} \right\rfloor \\       &= \Psi(M_2) 
			+ \left\lfloor \frac{|V_{G\setminus M_1}|}{2}  \right\rfloor -1\\
			&\geq \Psi(M_1) 
			+ \left\lfloor \frac{|V_{G\setminus M_1}|}{2} \right\rfloor -1 \\
			&\geq \Psi(G).
		\end{align*}
		Therefore $\Psi(M_1) = \Psi(M_2)$ and 
		$\Psi(G) = \Psi(M_2) 
		+ \left\lfloor \frac{|V_{G\setminus M_2}|}{2} \right\rfloor$.
		
		For the converse, observe that we would have
		\begin{align*}
			\Psi(G) &= \Psi(M_2) 
			+ \left\lfloor \frac{|V_{G\setminus M_2}|}{2} \right\rfloor \\
			&= \Psi(M_1) 
			+ \left\lfloor \frac{|V_{G\setminus M_1}|}{2} \right\rfloor - 1
		\end{align*}
		so that $\mpd_G(|M_1|)\leq \Psi(G)-\Psi(M_1)<\left\lfloor \frac{|V_{G\setminus M_1}|}{2}\right\rfloor$, and therefore $G$ is not weakly critical.
		
		For the final statement, write $\xi = \left\lfloor \frac{|V_{G\setminus M_1}|}{2} \right\rfloor$. Observe that 
		\[ \Psi(G) = \Psi(M_1) + \xi -1
		\leq |M_1| + \xi -1 \]
		and
		\[ |V| = |M_1| + |V_{G\setminus M_1}| 
		\geq |M_1| + 2\xi.\]
		If one vertex of each color is selected from $G$, we see that leaves at least $\xi+1$ vertices which may be chosen as $C$.
	\end{proof}
	
	It is worth pointing out that an analogous argument establishes the following characterization of critical graphs.
	
	\begin{theorem}    \label{thm: alt def of not critical}
		Let $G$ be a graph. Then $G$ is not critical if and only if 
		there exists induced subgraphs $M_1 \subseteq M_2 \subseteq G$ so that
		\begin{enumerate}
			\item $0<|V_{M_2 \setminus M_1}| \leq 2$,
			\item $\Psi(M_1) = \Psi(M_2)$,
			\item $\Psi(G) = \Psi(M_2) 
			+ \left\lceil \frac{|V_{G\setminus M_2}|}{2} \right\rceil$.
		\end{enumerate}
		Moreover, if $\ell$ is a maximal pseudocomplete coloring of $G$, 
		there exists $C\subseteq V_G$ with $|C| = \left\lceil \frac{|V_{G\setminus M_2}|}{2} \right\rceil+|V_{M_2 \setminus M_1}|-1$ and
		$\ell(G) = \ell(G\setminus C)$.
	\end{theorem}

	As a consequence of Theorem \ref{thm: alt def of not weakly critical}, we are now able to prove the following, which is as close to a converse of Corollary \ref{thm: G join H psi additive when critical} as possible.
	\begin{theorem} \label{thm: psi add gives at least a weak crit on one}
		Let $G$ and $H$ be graphs. If
		\[ \Psi(G) + \Psi(H) = \Psi(G\,\nabla H),\]
		then at least one of $G$ or $H$ is weakly critical. Moreover, if one of $G$ or $H$ is weakly critical and has a coloring of form (1) as in Theorem \ref{thm: weakly crit lots of edges}, then the other is weakly critical and does not have such a coloring.
	\end{theorem}
	
	\begin{proof}
		Let $G$ and $H$ be graphs with $\psi(G\,\nabla H)=\Psi(G)+\Psi(H)$.
		
		Suppose that both $G$ and $H$ are not weakly critical.
		Using Theorem \ref{thm: alt def of not weakly critical}, choose induced
		subgraphs $M_1 \subseteq M_2 \subseteq G$ and $N_1 \subseteq N_2 \subseteq H$. After possibly relabeling, suppose 
		$\left\lfloor \frac{|V_{G\setminus M_1}|}{2} \right\rfloor 
		\leq 
		\left\lfloor \frac{|V_{H\setminus N_1}|}{2} \right\rfloor$. 
		
		Color $H$ with a pseudocomplete coloring with $\Psi(H)$ colors.
		Choose $C\subseteq H$ with $|C| = \left\lfloor \frac{|V_{G\setminus M_1}|}{2} \right\rfloor$ so that all $\Psi(H)$ colors are still represented in $H\setminus C$. 
		
		Write $V_{G\setminus M_1} = P_0 \amalg P_1$ with $|P_0|=\left\lfloor \frac{|V_{G\setminus M_1}|}{2} \right\rfloor\leq|P_1|$. Color $V_{G\setminus M_1}$ with $\left\lfloor \frac{|V_{G\setminus M_1}|}{2} \right\rfloor$ new colors such that each color appears in both $P_0$ and $P_1$.
		
		Finally, swap the colors in $P_0$ with those in $C$ and color $M_1$ with a pseudocomplete coloring consisting of $\Psi(M_2)$ additional new colors. 
		
		By construction, this gives a pseudocomplete coloring with 
		\[   \Psi(H) + \Psi(M_2) + 
		\left\lfloor \frac{|V_{G\setminus M_2}|}{2} \right\rfloor + 1 \\
		=\Psi(H) + \Psi(G) + 1  \]
		colors. As this is a lower bound for $\Psi(G\,\nabla H)$, we get 
		$ \Psi(G) + \Psi(H) >  \Psi(G\,\nabla H)$
		which contradicts addivity of $\Psi$. 
		
		Thus one of $G$ or $H$ must be weakly critical. 
		
		Now, assume that $G$ is weakly critical and has a pseudocomplete coloring of type (1). Let $\ell$ be such a coloring, i.e. a coloring using $\Psi(G)=\left\lfloor\frac{\omega(G)+|V_G|}{2}\right\rfloor$ many colors with vertices $v_0$, $v_1$, and $v_2$ with $\ell(v_0)=\ell(v_1)=\ell(v_2)$. Note that $\ell$ uses all $\Psi(G)$ colors on $G\setminus\{v_0,v_1\}$.
		
		Assume that $H$ is either not weakly critical or weakly critical with a coloring of type (1). In either case, we will show that there are distinct $c_0,c_1\in V_H$ and a pseudocomplete coloring $\ell'$ of $H$ so that $\Psi(H)$ many colors appear in the complement of $\{c_0,c_1\}$. 
		
		If $H$ is not weakly critical, by Theorem \ref{thm: alt def of not weakly critical}, there is a pseudocomplete coloring $\ell'$ of $H$ using $\Psi(H)$ many colors and a subset $C$ of $V_H$ with $|C|=\left\lfloor\frac{V_{H\setminus N_1}}{2}\right\rfloor+1\geq 2$ such that each of the $\Psi(H)$ colors appears in the complement of $C$. Let $\{c_0,c_1\}\subseteq C$.
		
		If $H$ is weakly critical with a coloring of type (1), take $\ell'$ to be one such coloring and let $c_0$, $c_1$, and $c_2$ be the three vertices that share a color. Clearly $\ell'$ uses all $\Psi(H)$ colors on $H\setminus\{c_0,c_1\}$.  
		
		Now, color $G\,\nabla H$ as follows:  color $G\setminus\{v_0,v_1\}$ according to $\ell$; color $H\setminus\{c_0,c_1\}$ according to $\ell'$; color $c_0$ with $\ell(v_0)$; $v_0$ with $\ell'(c_0)$; and use one additional color to color both $v_1$ and $c_1$.  By construction, this is a pseudocomplete coloring consisting of $\Psi(G)+\Psi(H)+1>\Psi(G)+\Psi(H)$ many colors, establishing that $\Psi(G\,\nabla H)>\Psi(G)+\Psi(H)$, contradicting additivity of $\Psi$.
		Thus if $G$ is weakly critical and has a pseudocomplete coloring of type (1), then $H$ is weakly critical and does not have a pseudocomplete coloring of type (1).
	\end{proof}

\section{Remarks on $\nabla^k G$} \label{sec: nabla^k G}

If we restrict to the case that $G=H$, we can get a partial converse to Theorem \ref{thm: G join H psi additive when critical} with a stronger conclusion than that of Theorem \ref{thm: psi add gives at least a weak crit on one}. 

\begin{theorem} \label{thm: crit & psi additive for H=G}
	Let $G$ be a graph. Then $G$ is critical if and only if
	\[ \Psi(G \,\nabla G) = 2\Psi(G). \]
\end{theorem}

\begin{proof}
	Theorem \ref{thm: criteria for Psi to be additive over the join via mpd} shows that 
	$\Psi(G \,\nabla G) = 2\Psi(G)$ if and only if for all $k\leq |G|$, $\mpd_G(k) + \mpd_G(k)\geq k$. Since $\mpd_G(k)$ is an integer for all $k$, this happens if and only if $\mpd_G(k)\geq \left\lceil \frac{k}{2} \right\rceil$ for all $k\leq|G|$, \emph{i.e.}, $G$ is critical.
\end{proof}

\begin{remark} \label{rmk: false: converse of crit join crit is crit with G equals H}
	As noted in Remark \ref{rmk: false: converse of crit join crit is crit}, the converse of Theorem \ref{thm: crit join crit is crit} fails. Indeed the converse fails even under the assumption that $G=H$. Theorem \ref{thm: upper lower bounds for general psi of join} immediately implies that
	\begin{equation} \label{eqn: Psi of GjoinG}
		\Psi(G\,\nabla G) = \omega(G) + |V_G|
	\end{equation}
	for all $G$, so by Lemma \ref{lem: psi vs clique number inequality}, $G\,\nabla G$ is always critical.
\end{remark}
In fact, more is true.
\begin{theorem} \label{thm: k-join of G critical when it can be}
	Let $k\in\Z$ with $k\geq 2$. Then $\nabla^k G$ is critical if and only if 
	\[ k\,(\omega(G) + |V|) \]
	is even.
\end{theorem}

\begin{proof}
	If $k$ is even, apply Equation \ref{eqn: Psi of GjoinG} with $G$ replaced by $\nabla^{k/2} G$ and then use Lemma \ref{lem: psi vs clique number inequality} to get criticality.
	
	If $\omega(G) + |V|$ is even with $k\geq 3$, begin with a maximal clique, $K$, of $G$ and let $X=G\setminus K$. As $\omega(G) + |V|$ is even, so is $|X|=|V|-\omega(G)$. Divide $X$ into sets $X_1$ and $X_2$ of equal order, $q = \frac{|V|-\omega(G)}{2}$.
	
	In $\nabla^k G$, color $\nabla^k K$ with $k\omega(G)$ distinct colors. Next color $\nabla^k X_1$ with an additional $kq$ colors. Finally, color $\nabla^k X_2$ with a shift of the colors used in $\nabla^k X_1$. More precisely, color the $(i+1)$st copy of $X_2$ in $\nabla^k X_2$, $i$ viewed as an element of $\Z_k$, with the same colors used on the $i$th copy of $X_1$ in $\nabla^k X_1$. This ensures that the colors used in each copy of $X$ in $\nabla^k X$ are used in two different copies of $X$. It is straightforward to verify that this results in a pseudocomplete coloring of $\nabla^k G$. As
	\[ k\left(\omega(G) + \frac{|V|-\omega(G)}{2}\right) = k \frac{\omega(G) + |V|}{2}, \]
	Lemma \ref{lem: psi vs clique number inequality} again gives criticality.
	
	Finally, if $k(\omega(G) + |V|)$ is odd, Lemma \ref{lem: psi vs clique number inequality} shows that 
	$\nabla^k G$ is not critical.
\end{proof}

Note that, as criticality of $G$ requires $\omega(G) + |V|$ to be even, Theorem \ref{thm: upper lower bounds for general psi of join} says that $\nabla^k G$ is critical whenever parity makes criticality possible. 
See Theorem \ref{thm: k-join of G weakly critical for the rest} below when parity makes criticality impossible.

Theorem \ref{thm: k-join of G critical when it can be} allows us to generalize Theorem \ref{thm: crit & psi additive for H=G}. See Theorem \ref{thm: w crit if kpsi = psi nabla k with parity odd condition} for the analogue when $k\,(\omega(G) + |V|)$ is odd. 
Note that this, along with Theorem \ref{thm: w crit if kpsi = psi nabla k with parity odd condition} below, gives a correct proof of \cite{BRY2003} Corollary 11.
\begin{theorem} \label{thm: crit if kpsi = psi nabla k with parity even condition}
	Let $k\in\Z$ with $k\geq 2$ and $k\,(\omega(G) + |V|)$ even.
	Then $G$ is critical if and only if 
	\[ k\Psi(G) = \Psi(\nabla^k G). \]
\end{theorem}

\begin{proof}
	As Theorem \ref{thm: k-join of G critical when it can be} shows that $\nabla^k G$ is critical, it follows from Lemma \ref{lem: crit if 2psi=o+n} that $\Psi(\,\nabla^k G) = k\frac{\omega(G) + |V|}{2}$. Therefore, 
	$k\Psi(G) = \Psi(\nabla^k G)$ if and only if $\Psi(G) = \frac{\omega(G) + |V|}{2}$ if and only if $G$ is critical.
\end{proof}

Recall that Theorem \ref{thm: k-join of G critical when it can be} showed that, for $k\geq 2$, $\nabla^k G$ is critical if and only if $k(\omega(G) + |V|)$ is even. We now show that weakly critical fills in the rest of the range.
\begin{theorem} \label{thm: k-join of G weakly critical for the rest}
	Let $k\in\Z$ with $k\geq 2$. Then $\nabla^k G$ is weakly critical.
\end{theorem}

\begin{proof}
	By Theorem \ref{thm: k-join of G critical when it can be} and the comment following Theorem \ref{thm: weak critical numerical equivalence}, it remains to show that $\nabla^k G$ is weakly critical when $k(\omega(G) + |V|)$ is odd with $k\geq 3$. 
	
	Begin with a maximal clique, $K$, of $G$ and let $X=G\setminus K$. As $\omega(G) + |V|$ is odd, so is $|X|=|V|-\omega(G)$. Divide $X$ into sets $X_1$ and $X_2$ of equal order, $q = \frac{|V|-\omega(G)-1}{2}$, and a singleton vertex, $v_0$. 
	
	In $\nabla^k G$, as in Theorem \ref{thm: k-join of G critical when it can be}, color $\nabla^k K$ with $k\omega(G)$ distinct colors, $\nabla^k X_1$ with an additional $kq$ colors, and $\nabla^k X_2$ with a shifted coloring of $\nabla^k X_1$. Finally color the copies of $v_0$ in $\nabla^k G$ with an additional 
	$\frac{k-1}{2} $ colors. 
	It is straightforward to verify that this generates a pseudocomplete coloring of $\nabla^k G$. As
	\begin{align*}
		k\left(\omega(G) + \frac{|V|-\omega(G)-1}{2}\right) + \frac{k-1}{2} &= k \frac{\omega(G) + |V|}{2} -\frac{1}{2} \\
		&= \left\lfloor \frac{k(\omega(G) + |V|)}{2} \right\rfloor ,
	\end{align*}
	Theorem \ref{thm: weak critical numerical equivalence} gives weak criticality. 
\end{proof}

Theorem \ref{thm: k-join of G weakly critical for the rest} now allows us to address the analogue of Theorem \ref{thm: crit if kpsi = psi nabla k with parity even condition} when $k\,(\omega(G) + |V|)$ is odd.
\begin{theorem} \label{thm: w crit if kpsi = psi nabla k with parity odd condition}
	Let $k\in\Z$ with $k\geq 3$ and $k\,(\omega(G) + |V|)$ odd.
	Then $G$ is not critical. However, $G$ is weakly critical if and only if 
	\[ k\Psi(G) +\left\lfloor \frac{k}{2} \right\rfloor = \Psi(\nabla^k G). \]
\end{theorem}

\begin{proof}
	The fact that $G$ is not critical follows by parity from Lemma \ref{lem: crit if 2psi=o+n}.
	As Theorem \ref{thm: k-join of G weakly critical for the rest} shows that $\nabla^k G$ is weakly critical, it follows from Theorem \ref{thm: weak critical numerical equivalence} that 
	$\Psi(\,\nabla^k G) = \left\lfloor \frac{k(\omega(G) + |V|)}{2} \right\rfloor$. 
	However, we know that $G$ is weakly critical iff and only if 
	$\Psi(G) = \left\lfloor \frac{\omega(G) + |V|}{2} \right\rfloor$. As this is equivalent to 
	\begin{align*}
		k\Psi(G) &= k \left\lfloor \frac{\omega(G) + |V|}{2} \right\rfloor \\
		&= \frac{k(\omega(G) + |V|)}{2} - \frac{k}{2} \\
		&= \left\lfloor \frac{k(\omega(G) + |V|)}{2} \right\rfloor
		- \left\lfloor \frac{k}{2} \right\rfloor,
	\end{align*}
	we are done.
\end{proof}

\section*{Acknowledgments}

The third author and the Kalasalingam Academy of Research and Education promote the Sustainable Development Goal of ensuring inclusive and equitable quality education and promoting lifelong learning opportunities for all.

In addition, the third author is thankful to the Department of 
Mathematics, Baylor University, for the support given during his visit.

\bibliographystyle{abbrvnat}
\bibliography{refs.bib}
\end{document}